\makeatletter
\IfFileExists{./numapde-latex/numapde-packages.sty}{\providecommand*{\input@path}{}\edef\input@path{{./numapde-latex/}\input@path}}{}
\makeatother

\documentclass[%
]{article}

\usepackage{microtype}

\usepackage[%
 accepted
]{icml2021}

\usepackage{etoolbox}
\makeatletter
\patchcmd{\icmltitle}{\bf}{\bfseries}{}{}
\patchcmd{\icmltitle}{\bf}{\bfseries}{}{}
\patchcmd{\icmltitle}{\bf}{\bfseries}{}{}
\patchcmd{\icmltitle}{\bf}{\bfseries}{}{}
\patchcmd{\icmlauthor}{\bf}{\bfseries}{}{}
\patchcmd{\icmlauthor}{\bf}{\bfseries}{}{}
\patchcmd{\icmlaffiliation}{\bf}{\bfseries}{}{}
\patchcmd{\icmladdress}{\bf}{\bfseries}{}{}
\patchcmd{\printAffiliationsAndNotice}{\bf}{\bfseries}{}{}
\patchcmd{\printAffiliationsAndNotice}{\bf}{\bfseries}{}{}
\patchcmd{\abstract}{\bf}{\bfseries}{}{}
\patchcmd{\section}{\bf}{\bfseries}{}{}
\patchcmd{\subsection}{\bf}{\bfseries}{}{}
\patchcmd{\paragraph}{\bf}{\bfseries}{}{}
\patchcmd{\subparagraph}{\bf}{\bfseries}{}{}
\makeatother

\RequirePackage[algo2e,ruled,vlined,linesnumbered]{algorithm2e}
\usepackage{numapde-semantic}
\usepackage{numapde-style}
\usepackage{numapde-local}

\IfFileExists{numapde-uncolorizeHyperrefIfChanges.sty}{\RequirePackage{numapde-uncolorizeHyperrefIfChanges}}{}

\begin{document}

\twocolumn[

\icmltitle{Optimal Prediction using Learning and Shape Optimization}
\icmltitlerunning{}

\icmlsetsymbol{equal}{*}

\begin{icmlauthorlist}
\icmlauthor{M. Sajjad Edalatzadeh}{1} 
\icmlauthor{Roland Herzog}{1} 
\end{icmlauthorlist}
\icmlaffiliation{1}{Technische Universität Chemnitz, Faculty of Mathematics, 09107 Chemnitz, Germany}
\icmlcorrespondingauthor{M. Sajjad Edalatzadeh}{msedalatzadeh@gmail.com}

\icmlkeywords{learning, prediction, shape optimization, distributed parameter systems, sensors}

\vskip 0.3in
]

\printAffiliationsAndNotice{} 

\begin{abstract}
This paper investigates the problem of optimal predictor design for distributed parameter systems using neural networks and shape optimization. 
Sensors with various shapes are placed on the domain of the distributed parameter system. 
Data provided by these sensors are fed into a re-constructor to generate a full state of the system. 
After that, a trained neural-network predictor produces a prediction of the state at future time steps. 
The cost of prediction is defined as the weighted sensor area plus the squared norm of the prediction error.
The location and shape of a sensor influences the prediction cost as well as the predictor performance.
With the aid of the gradient of the network with respect to its inputs, an outer optimization layer is augmented to find optimized locations and shapes of the sensors. 
Simulation results show good agreement between the predicted and the true states as well as a significant reduction in cost by sensors with optimized locations and shapes.
\end{abstract}


\section{Introduction}
\label{section:introduction}

Estimation is of fundamental importance in many disciplines. 
For example, it is important in climate science to predict the temperature changes in the Earth atmosphere \cite{ParkKimLeeKimSongKim:2019:1}. 
Estimation is also used to find sea surface temperature in order to predict weather \cite{NielsenEnglystHoyerPedersenGentemannAlerskansBlockDonlon:2018:1}. 
In economics, estimation is used to predict investment risks. 
In control theory, estimation is used to feedback from the state of the system at locations not accessible by sensors. 
Most commonly, estimation is used to reconstruct the state of a system modeled by partial differential equations such as the heat and wave equations.

In recent years, machine learning methods have been used to extract predictive models from large sets of data provided by sensors. 
These methods can further find a symbolic representation of the model (symbolic regression) and also discover the hidden physics. 
Among the machine learning methods, deep neural networks have gained considerable attention \cite{LeCunBengioHinton:2015:1}. 
Evolutionary algorithms \cite{VaddireddyRasheedStaplesSan:2020:1}, compressive sensing \cite{CandesWakin:2008:1}, and sparse
optimization \cite{CandesWakinBoyd:2008:1} are also machine learning algorithms used for logistic regression. 
In \cite{LongLuMaDong:2018:1}, a PDE-net approach is proposed to predict the dynamics of systems modeled by partial differential equations (PDEs) from sensor data. 
The network is also able to uncover the underlying PDE model.
In \cite{LuZhongLiDong:2018:1}, it is shown that deep neural networks such as ResNet, PolyNet, FractalNet, and RevNet can be interpreted as different numerical discretizations of differential equations. 
This interpretation is further used to design new deep neural networks. 
In \cite{LiuLinZhangSu:2010:1,LiuLinZhangTangSu:2013:1}, training data for image restoration is used to create a learning-based PDE model for computer vision tasks. 
The problem is stated as an optimal control problem where the inputs and outputs are training images. 
The authors show the effectiveness of their model by numerical experiments on image denoising and inpainting problems in image restoration. 

The shape optimization of sensors and actuators in the context of PDEs has been discussed in relatively few works. 
In \cite{PrivatTrelatZuazua:2013:1}, the optimal shape and position of an actuator for the wave equation in one spatial dimension are discussed. 
An actuator is placed on a subset $\omega\in [0,\pi]$ with a constant Lebesgue measure $L\pi$ for some $L\in (0,1)$. 
The optimal actuator minimizes the norm of a Hilbert Uniqueness Method (HUM)-based control; such control steers the system from a given initial state to zero state in finite time. 
In \cite{PrivatTrelatZuazua:2017:1}, the optimal actuator shape and position for linear parabolic systems are discussed. 
This paper adopts the same approach as in \cite{PrivatTrelatZuazua:2013:1} but with initial conditions that have randomized Fourier coefficients. 
The cost is defined as the average of the norm of HUM-based controls. 
In \cite{KaliseKunischSturm:2018:1}, optimal actuator design for linear diffusion equations has been discussed. 
A quadratic cost function is considered, and shape and topological derivatives of this function are derived. 
Numerical results show significant improvement in the cost and performance of the control. 
In \cite{EdalatzadehKaliseMorrisSturm:2019:1}, the optimal shape of actuators for vibration control of flexible beams is investigated using a linear-quadratic performance index and shape optimization.  

Optimal sensor design problems are in many ways similar to the optimal actuator design problem. 
In \cite{PrivatTrelatZuazua:2014:1}, optimal sensor shape design has been studied where the observability is being maximized over all admissible sensor shapes. 
Optimal actuator design problems for nonlinear distributed parameter systems have also been studied. 
In \cite{EdalatzadehMorris:2019:1}, it is shown that under certain conditions on the nonlinearity and the cost function, an optimal input and actuator design exist, and optimality conditions are derived. 
Results are applied to the nonlinear railway track model as well as to the semi-linear wave model in two spatial dimensions. 
Numerical techniques to calculate the optimal actuator shape design are mostly limited to linear quadratic regulator problems, see, \eg, \cite{AllaireMuenchPeriago:2010:1}. 
For controllability-based approaches, numerical schemes have been studied in \cite{MuenchPeriago:2011:1,Muench:2007:1,MuenchPeriago:2013:1}.

In this paper, we combine shape optimization with learning-based PDE models to find the optimal shape of sensors, as well as an estimator of the full state of the PDE. 
In \cref{section:estimator_design}, the neural-network predictor is formulated, which recovers the full PDE state from limited sensor data and predicts the state at the subsequent time step. 
\Cref{section:shape_optimization} discusses the shape optimization of sensors. 
Simulations results are presented in \cref{section:simulation_results}. 
Concluding remarks are given in \cref{section:conclusion}.

\section{Estimator design}
\label{section:estimator_design}
Consider a first-order dynamical system with state $z(t)$, initial condition $z_0$, state space $X$, and linear or nonlinear operator $F$ on $X$:
\begin{equation}
\label{eq:first-order_system}
\paren[auto]\{.{%
	\begin{aligned}
	\dot{z}(t) & = F(z)
	, 
	\\
	z(0) & = z_0
	.
	\end{aligned}
}
\end{equation} 
In our applications, $F$ is a second-order partial differential (PDE) operator acting on functions in one or two spatial dimensions.
The evolution of the state $z(t)$ in discrete time with sufficiently small time increment $\Delta t$ can be described, using for simplicity an explicit Euler method, as
\begin{equation}
\label{eq:first-order_system_discrete_in_time}
z_{k+1} = z_k + F(z_k) \Delta t
\end{equation}
with initial value $z_0$.
For numerical computation, we also need to discretize in space.
For simplicity of notation, we currently assume $X$ to be a space of functions on a one-dimensional interval, which we discretize with a mesh of $c$ equidistant points.
The results provided in \cref{section:simulation_results} include two-dimensional examples as well.
Discretizing $F$ accordingly, here using a finite difference method, we obtain from \eqref{eq:first-order_system_discrete_in_time} the fully discrete evolution
\begin{equation}
\label{eq:first-order_system_fully_discrete}
\tz_{k+1} = \tz_k + \tF(\tz_k) \Delta t
.
\end{equation}
Each element in the vector $\tz_k$ represents an approximation of the true state $z$ at one of the spatial grid points and time $t = k \, \Delta t$.

A neural-network predictor takes the current state vector $\tz_k$ (or a reconstruction thereof) and predicts the solution at the next time step $\tz_{k+1}$. 
The estimator network is trained over simulation data or acquired real-life data. 
In most cases, the state $\tz_k$ cannot be fully measured by an array of sensors. 
Therefore, the data provided by the sensors available represents an incomplete state vector and the full state vector is reconstructed first before being fed into the predictor network.

\begin{figure}[H]
	\caption{Illustration of the reconstruction function $f_r$.}
	\centering
	\includegraphics[width=0.5\textwidth]{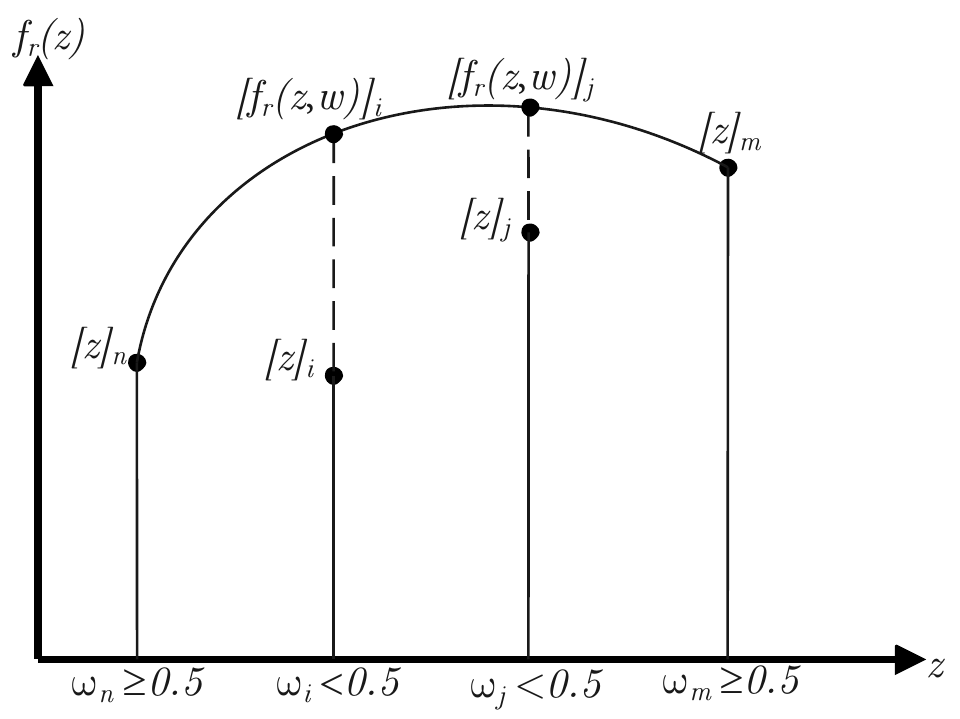}
	\label{fig-recons}
\end{figure}

The reconstruction of an incomplete state reading works as follows.
A sensor arrangement is encoded in a $c$-dimensional vector $\omega$. 
Each element of $\omega$ indicates the weight of the presence of a sensor at the respective vertex.
Weights less than 0.5 will be treated as no sensor present.
The function $f_r \colon \R^c\times \R^c \to \R^c$ provides an approximate reconstruction based on an interpolation of neighboring sensors measurements, \Cref{fig-recons}. 
That is, we have $\tz \approx f_r(\tz,\omega)$.
A cubic spline interpolation is used in the one-dimensional simulations. 
This interpolation is a suitable choice for PDEs with second-order spatial derivatives. 



The function $f_n \colon \R^c \to \R^c$ provides a neural-network prediction $\tZ_{k+1}$ of the state $\tz_{k+1}$ using reconstructed data; that is
\begin{equation}
\label{eq:reconstruction_followed_by_prediction}
\tz_{k+1} 
\approx 
\tZ_{k+1} 
\coloneqq
f_n(f_r(\tz_k,\omega))
.
\end{equation}
Different types of neural networks lead to different formulations for the predictor function $f_n$. 
Examples are given in \cref{section:simulation_results}.

\section{Sensor shape optimization}
\label{section:shape_optimization}

In this section, we consider an optimization of the sensor configuration encoded in the weight vector~$\omega$.
Clearly, from the prediction point of view, it is beneficial to have sensors everywhere.
In this case, we would have $\tz_k = f_r(\tz_k,\omega)$ since no interpolation is necessary.
On the other hand, we associate a cost with each sensor placed.
Our goal is to obtain a compromise between the number of sensors and the accuracy of the predictions $\tZ_{k+1}$ of the state $\tz_{k+1}$, based on the partially observed previous state $\tz_k$.
We express this goal in terms of the cost function
\begin{multline}
	\label{eq:cost_function}
	J(\omega)
	=
	\alpha \, \sum_{i=1}^c[\omega]_i
	\\
	+
	\sum_{k=0}^{K-1} \sum_{i=1}^c \paren[big][]{{f_n(f_r(\tz_k,\omega)) - \tz_{k+1}}}_i^2 \, h_i \, \Delta t
	.
\end{multline}
Here $K$ denotes the number of grid points in time. 
Notice that the second term in \eqref{eq:cost_function} represents a discretized version of the $L^2$ mean squared error.
The sequence $h_i$ depends on the spatial mesh size and the quadrate scheme used; for instance, $h_i = h/2,h,\cdots,h,h/2$ implements the trapezoidal rule with mesh size $h$. 
A similar formula holds in two spatial dimensions.

The derivative of $J(\omega)$ with respect to (\wrt) the sensor weight vector $\omega$ is
\begin{multline}
	\label{eq:derivative_of_cost_function}
	[J'(\omega)]_i
	=
	\alpha 
	\\
	+
	2\sum_{k=0}^{K-1} \sum_{i=1}^c \paren[big][]{\paren[big](){f_n(f_r(\tz_k,\omega)) - \tz_{k+1}}^\transp \tZ_{k+1}'}_i \, h_i \, \Delta t
	.
\end{multline}
Here $\tZ_{k+1}'$ is the derivative of $\tZ_{k+1} = f_n(f_r(\tz_k,\omega))$ \wrt\ $\omega$.
From the chain rule, we obtain
\begin{equation}
\label{eq:chain_rule}
\tZ_{k+1}'
=
f_n'(f_r(\tz_k,\omega)) \, f_r'(\tz_k,\omega)
,
\end{equation}
where the right hand side is a product of two $c \times c$ matrices.
The first factor represents the derivative of the estimator network's output \wrt\ its input vector, while the second factor is the derivative of the outcome of the spline interpolation \wrt\ $\omega$.
The latter derivative is zero except when there are components of $\omega$ right at the threshold $\omega = 0.5$.
In this case, the derivative becomes a Dirac function.
Simulation results, however, indicate that the first term in \eqref{eq:chain_rule} is sufficient to provide search directions suitable to reduce the cost function.
Therefore, we use the approximation
\begin{equation}
\label{eq:approximate_}
\tZ_{k+1}'
\approx
f_n'(f_r(\tz_k,\omega)) 
\end{equation}
to evaluate the derivative $J'(\omega)$ in \eqref{eq:derivative_of_cost_function}.
Its transpose, the approximate gradient $\nabla J(\omega)$, is then used in a trust-region optimization algorithm to find an improved sensor configuration vector $\omega$.

\begin{algorithm2e}
	\caption{Optimal learning-based predictor design for distributed parameter systems}
	\label{algorithm:deep_learning_and_shape_optimization}
	\SetAlgoLined
	Generate training data {from the simulated }solution $\tz$\\
	Build and train an estimator neural network model $f_n$\\
	Randomly select an initial weight vector $\omega_0$\\
	\For{$j=0, \ldots, J$}{
		\For{$k=0, \ldots, K$}{
			Evaluate $\tZ_{k+1} = f_n(f_r(\tz_k,\omega_j))$\\
			Evaluate $\tZ'_{k+1} \approx  f'_n(f_r( \tz_k, \omega_j))$}
		Find $\omega_{j+1}{ \in [0,1]^c}$ using a constrained optimization method{, which utilizes evaluations of }$J(\omega_j)$ and $J'(\omega_j)$ based on \eqref{eq:cost_function} and \eqref{eq:derivative_of_cost_function}\\
		$j\gets j+1$
	}
\end{algorithm2e}

\section{Simulation results}
\label{section:simulation_results}

In this section, we provide numerical results for a range of different PDE models.
Our goal is to demonstrate that, in each case, optimized sensor configurations and properly trained predictor networks $f_n$ can obtain sufficiently accurate predictions of the state vector.
To this end, we show that the cost $J(\omega)$ for an optimized weight $\omega$ is significantly lower than the cost for the all-sensor case $[\omega]_i = 1$, $i=1,2,\ldots,c$, for which the lower bound $J(\omega) \ge \alpha \, c$ is valid.
The overall procedure is described in \cref{algorithm:deep_learning_and_shape_optimization}.

For each example, we build an estimator neural network $f_n$ using the deep learning package Keras (version~2.4.0). 
We employ a linear activation function and the mean squared error (MSE) loss function in the networks and use the Adam optimizer for training.
The Adam optimizer and the MSE loss function yielded better performance compared to the rest of available optimization algorithms and loss functions. 
For PDE models in one space dimension, the network layout consists of two dense layers, where each layer has as many neurons as grid nodes. 
For PDE models in two space dimensions, a convolutional neural network with one layer, one filter and kernel size $3\times 3$ and zero padding is considered. 
In each case, the gradient of the network is calculated using the method \textit{backend} in Keras.

In the implementation of the reconstructor function~$f_r$, the sensor data is reconstructed using the Python package \textit{interpolate}. 
For PDE models in one space dimension, \textit{CubicSpline} is used, whereas \textit{interp2d} is chosen in two space dimensions.

Simulations are presented for a variety of PDE models in the following subsections. 
In figures and formulas, $u_r$ denotes the true (simulated) solution, which corresponds to $\tz$ in the notation of \cref{section:estimator_design,section:shape_optimization}.
Moreover, $u_p$ denotes the predicted solution corresponding to $\tZ$ according to \eqref{eq:reconstruction_followed_by_prediction}.

\subsection{One-dimensional heat equation}
\label{subsection:1Dlinearheat}
Consider the following heat equation in one space dimension over $x\in[0,1]$
\begin{equation}\label{1Dlinearheat}
\begin{cases}
u_t (x,t) = \kappa \, u_{xx}(x,t),\\
u(x,0) = u_0,\\
u(0,t) = 0,\\
u_x (1,t) = 0.
\end{cases}
\end{equation}

\begin{figure}[htp]
	\caption{Snapshots of prediction performance for the linear one-dimensional heat equation \eqref{1Dlinearheat} with initial condition $u_0(x) = x^2(x-1)^2(x-1/2)^2$, scaled to $[0,1]$. The solid blue line is the solution; the {orange} broken line is the prediction; green squares show the optimized sensor arrangement.}
	\label{fig-snapshots-heat-IC1}
	\centering
	\includegraphics[width=0.5\textwidth]{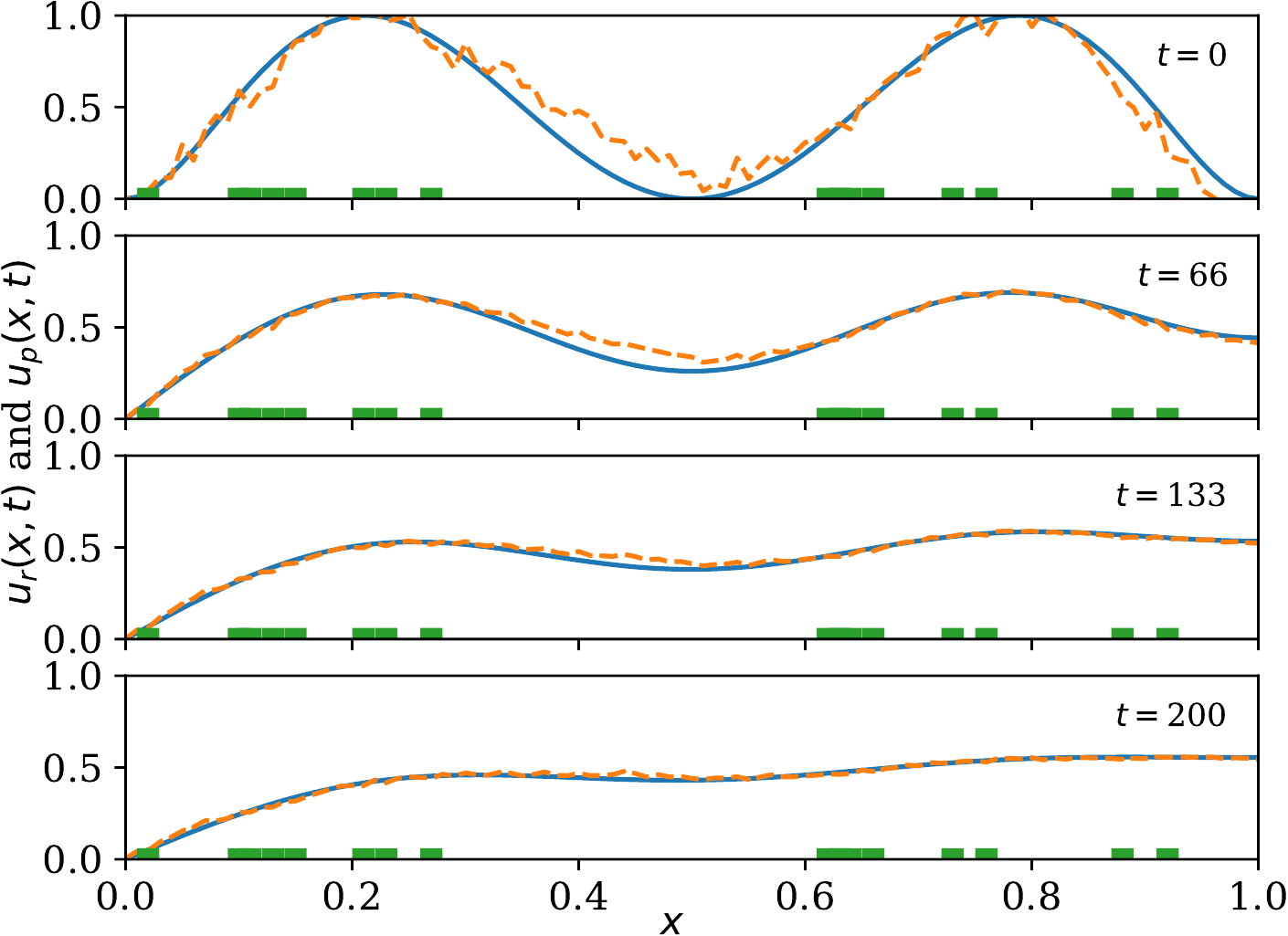}
\end{figure}

\begin{figure}[htp]
	\caption{Same as \cref{fig-snapshots-heat-IC1} but with initial condition $u_0(x)=x^2(x-1)^2(x-1/4)^2$, scaled to $[0,1]$.}
	\label{fig-snapshots-heat-IC2}
	\centering	\includegraphics[width=0.5\textwidth]{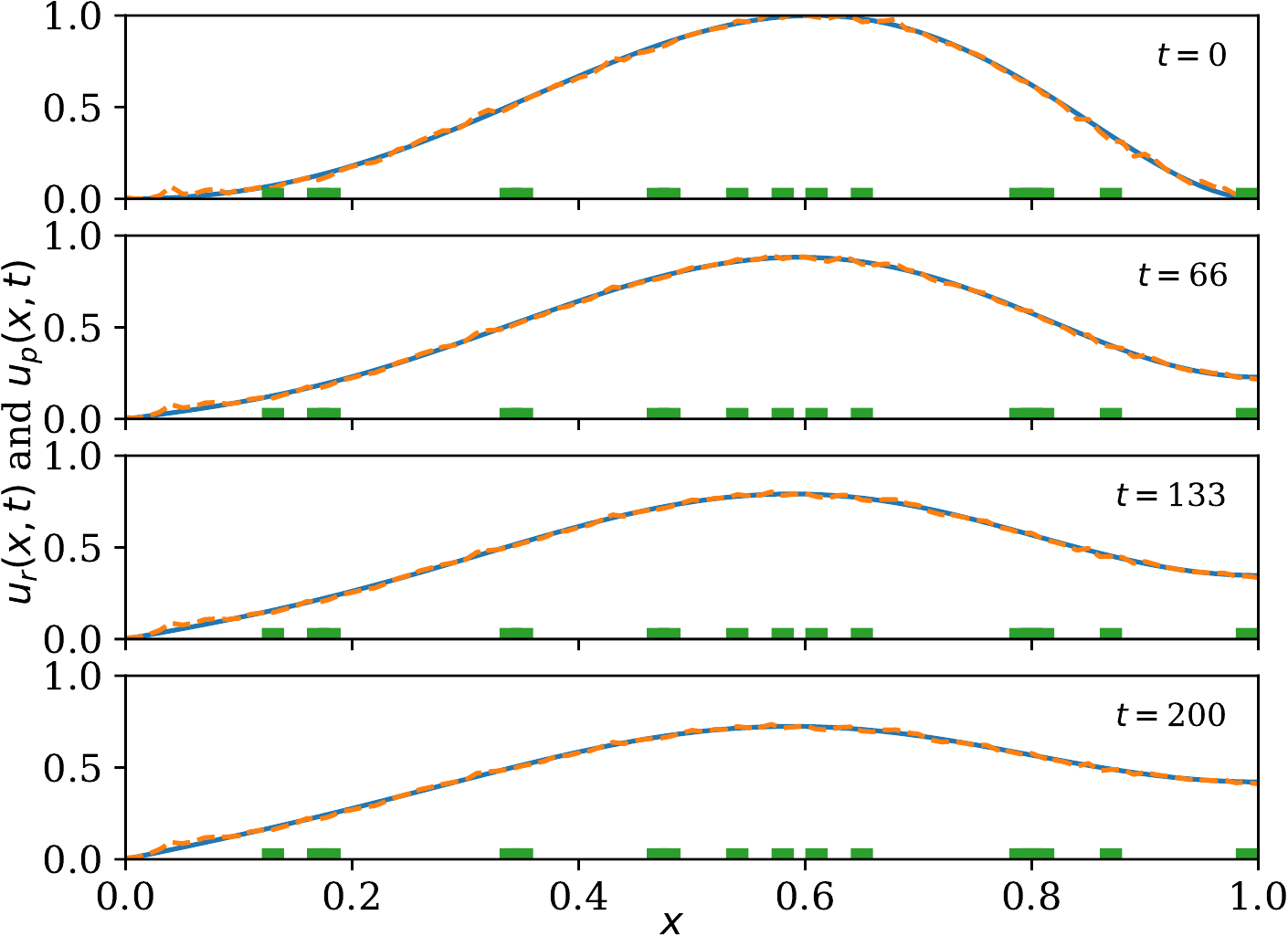}
\end{figure}

\begin{figure}[htp]
	\caption{Same as \cref{fig-snapshots-heat-IC1} but with initial condition $u_0(x)=x^2(x-1)^2(x+1/2)^2$, scaled to $[0,1]$.}
	\label{fig-snapshots-heat-IC3}
	\centering
	\includegraphics[width=0.5\textwidth]{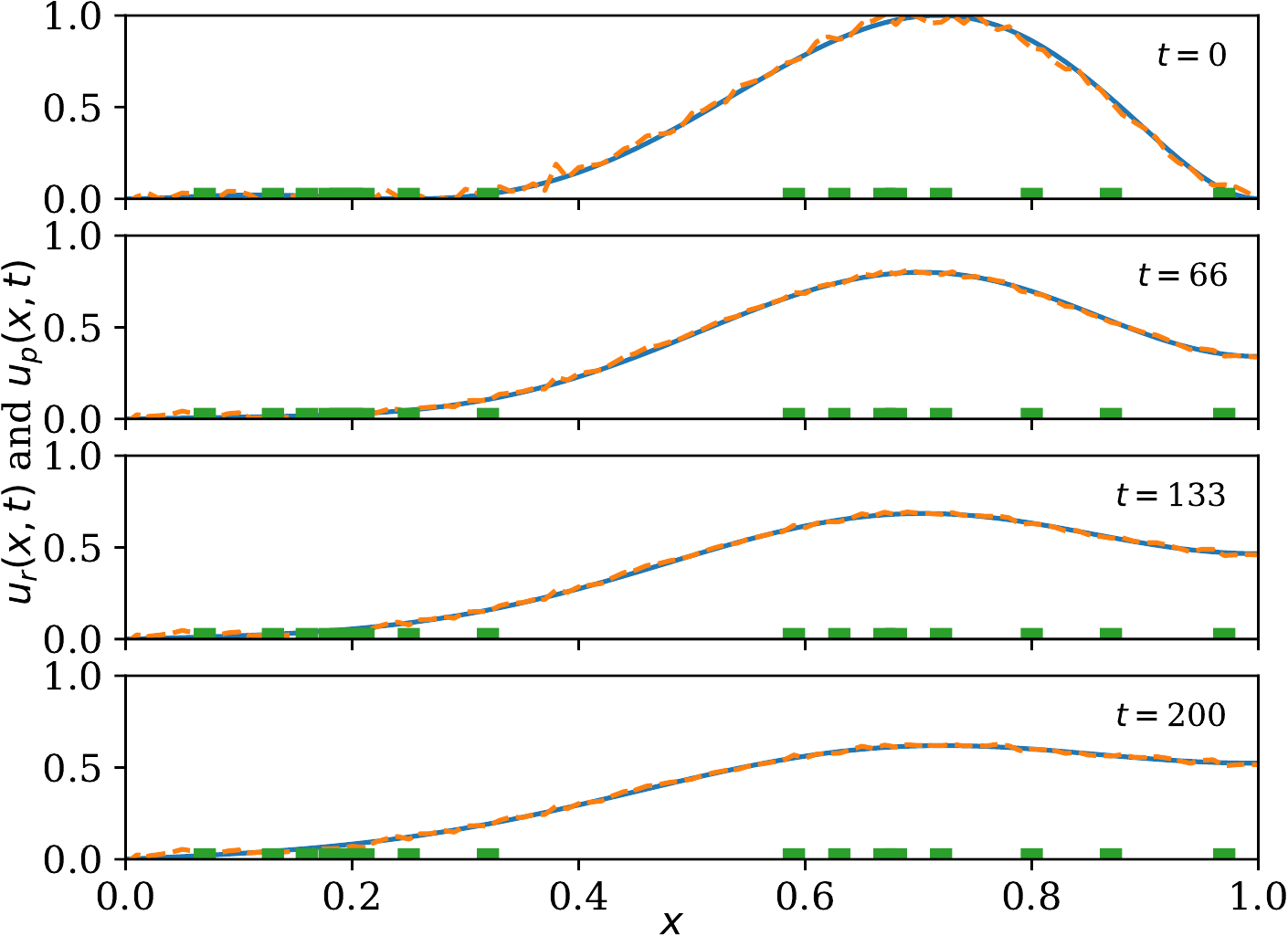}
\end{figure}

\begin{figure}[htp]
	\caption{Same as \cref{fig-snapshots-heat-IC1} but with initial condition $u_0(x)=H(x-1/2)$, where $H$ denotes Heaviside function, scaled to $[0,1]$.}
	\label{fig-snapshots-heat-step}
	\centering
	\includegraphics[width=0.5\textwidth]{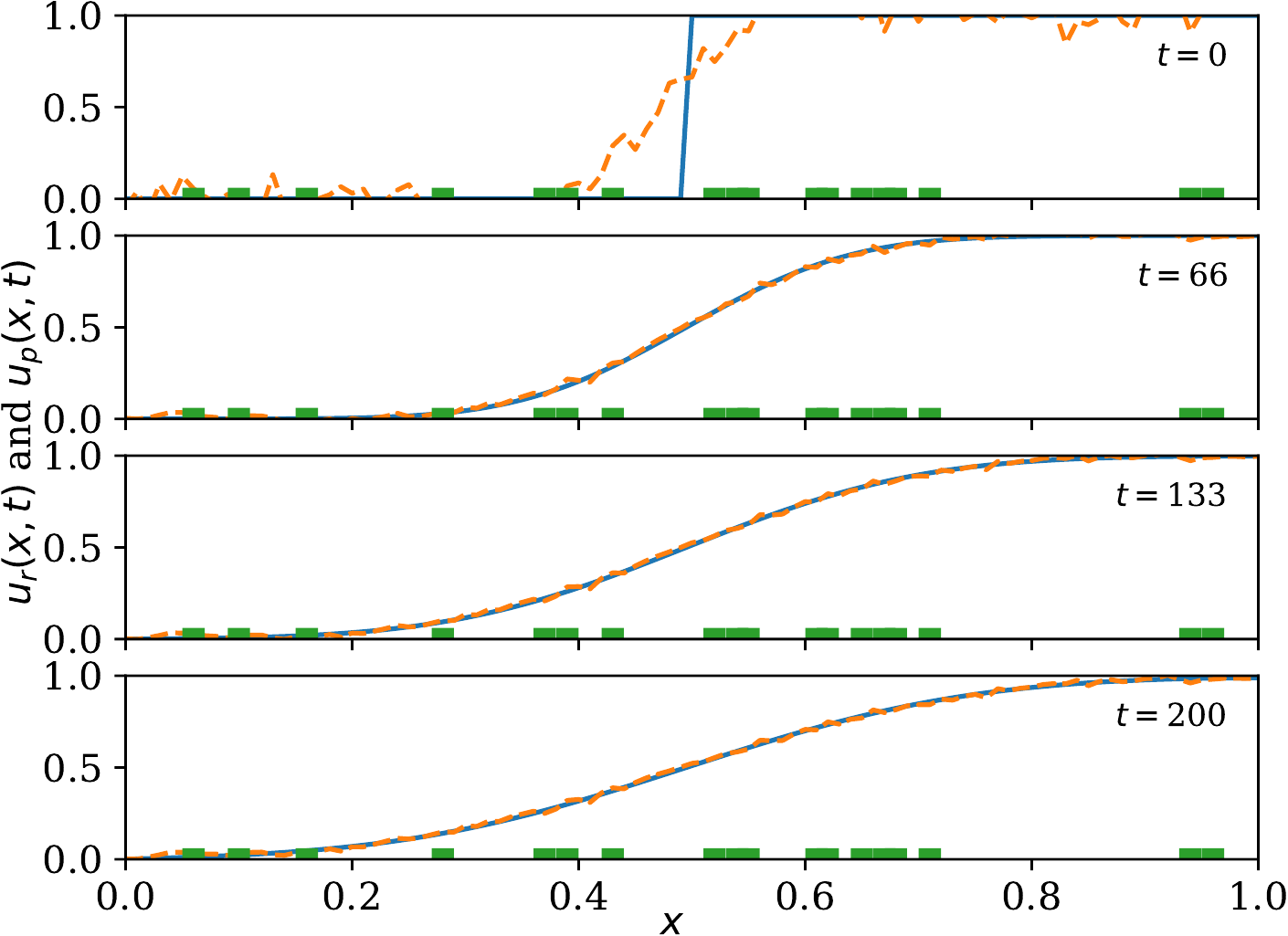}
\end{figure}

A forward-in-time and space-centered finite difference method is used to discretize the model and extract the solution data $\tz_k$. 
The predictor neural network is then trained on the solution data. 
The estimation performance is shown in \cref{fig-snapshots-heat-IC1,fig-snapshots-heat-IC2,fig-snapshots-heat-IC3,fig-snapshots-heat-step} for various initial conditions $u_0$. 
The sensor weight $\alpha$ is set to 5; mesh size $h$ is set to $10^{-2}$ (resulting in $c = 101$); time increment $\Delta t$ is set to $0.1$; and conductivity $\kappa$ is set to $10^{-4}$.

\begin{figure}[htp]
	\caption{$L^1$-norm of prediction error over time for the linear one-dimensional heat equation \eqref{1Dlinearheat}.}
	\label{fig-l1norm-linearheat}
	\centering
	\includegraphics[width=0.5\textwidth]{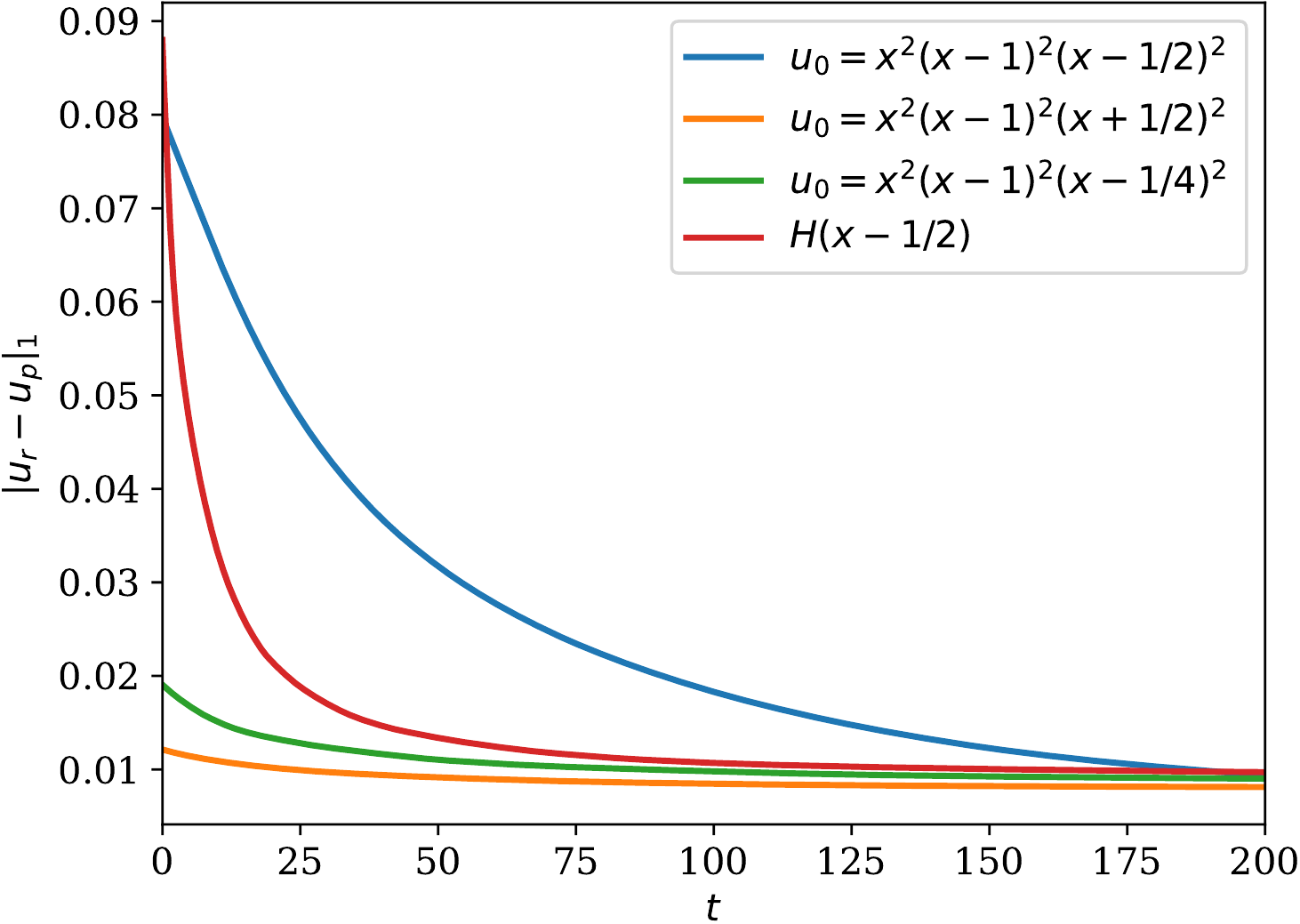}
\end{figure}

\cref{fig-l1norm-linearheat} shows the $L_1$-norm of the error over time, \ie,
\begin{equation}
\norm{u_r-u_p}_1
= 
\int_0^1 \abs{u_r(x,t)-u_p(x,t)} \d x
,
\end{equation} 
{which is small compared to the solution although only a fraction of the sensors is being used in each case, see \cref{fig-snapshots-heat-IC1,fig-snapshots-heat-IC2,fig-snapshots-heat-IC3,fig-snapshots-heat-step}.}
{In these figures, 17\% of the PDE region is covered with sensors.} 


\begin{figure}[htp]
	\caption{Reduction in cost for the linear one-dimensional heat equation \eqref{1Dlinearheat} and initial conditions given in \cref{fig-snapshots-heat-IC1,fig-snapshots-heat-IC2,fig-snapshots-heat-IC3,fig-snapshots-heat-step}. }
	\label{fig-cost-heat-linear}
	\centering
	\includegraphics[width=0.5\textwidth]{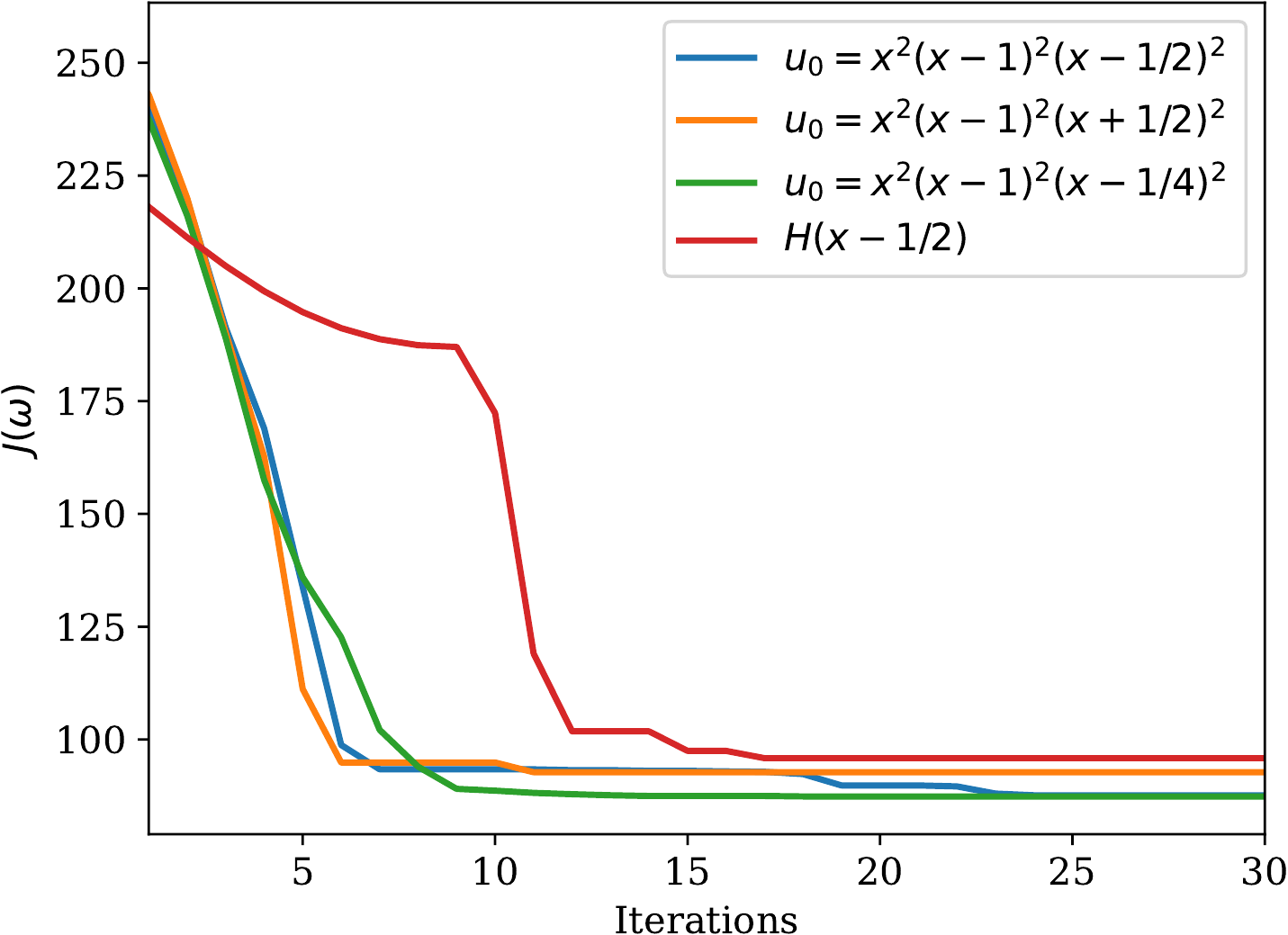}
\end{figure}

The reduction in cost over iteration is shown in \cref{fig-cost-heat-linear}.

\subsection{One-dimensional wave equation}
\label{subsection:1Dlinearwave}
Consider the following linear one-dimensional wave equation over $x\in [0,1]$
\begin{equation}\label{1Dlinearwave}
\begin{cases}
u_{tt} (x,t) = \lambda \, u_{xx}(x,t),\\
u(x,0) = u_0,\\
u_t(x,0) = 0,\\
u(0,t) = u(1,t) = 0.
\end{cases}
\end{equation}

\begin{figure}[h!]
	\caption{Snapshots of prediction performance for the linear one-dimensional wave equation \eqref{1Dlinearwave} with initial condition $u_0(x) = x^2(x-1)^2(x-1/2)^2$, scaled to $[0,1]$. The solid blue line is the solution; the {orange} broken line is the prediction, and green squares show the optimized sensor arrangement.}
	\label{fig-snapshots-wave-IC1}
	\centering
	\includegraphics[width=0.5\textwidth]{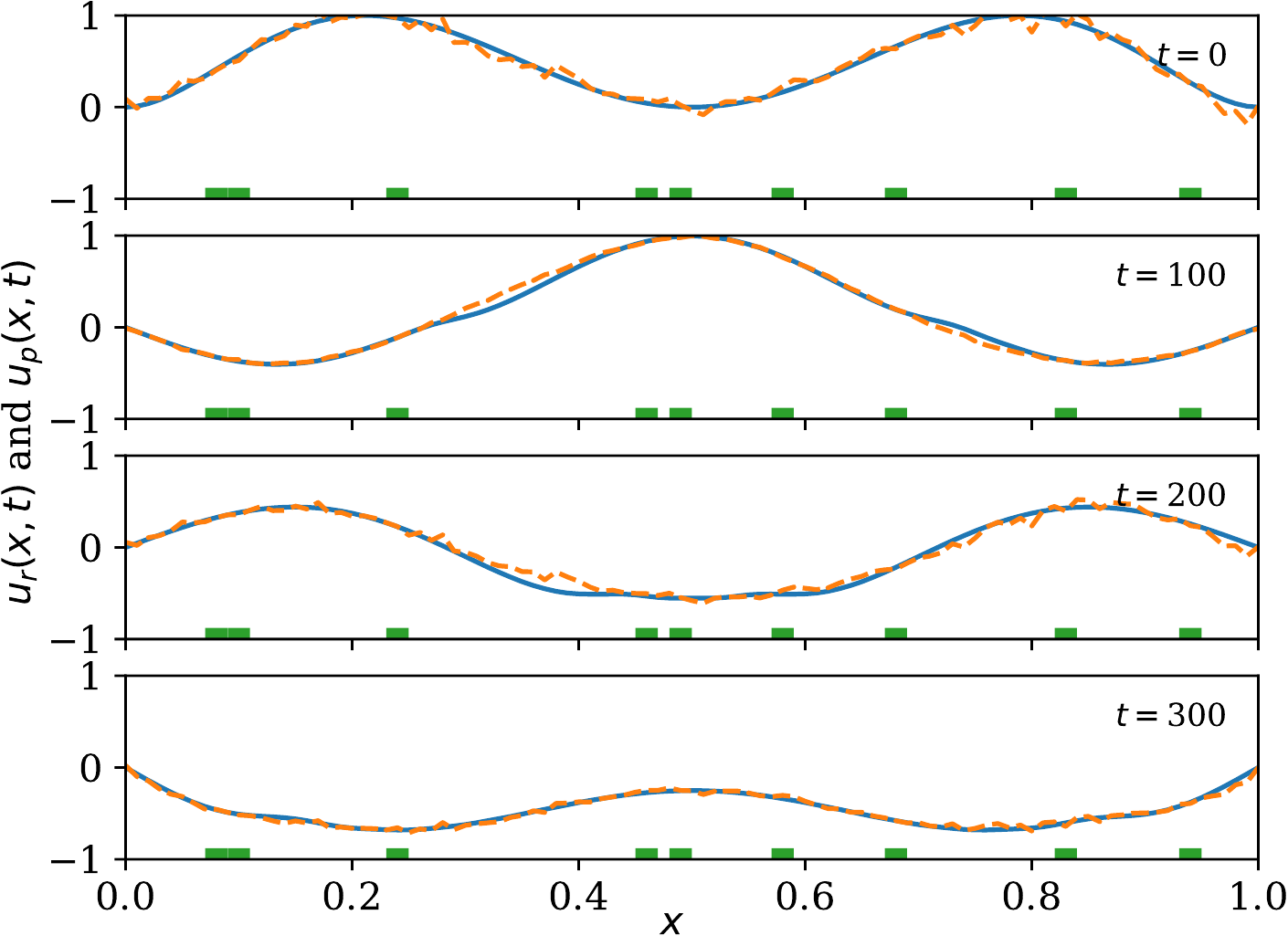}
\end{figure}

\begin{figure}[h!]
	\caption{Snapshots of prediction performance for the linear one-dimensional wave equation \eqref{1Dlinearwave} with initial condition $u_0(x) = x^2(x-1)^2(x+1/2)^2$, scaled to $[0,1]$. The solid blue line is the solution; the {orange} broken line is the prediction, and green squares show the optimized sensor arrangement.}
	\label{fig-snapshots-wave-IC2}
	\centering
	\includegraphics[width=0.5\textwidth]{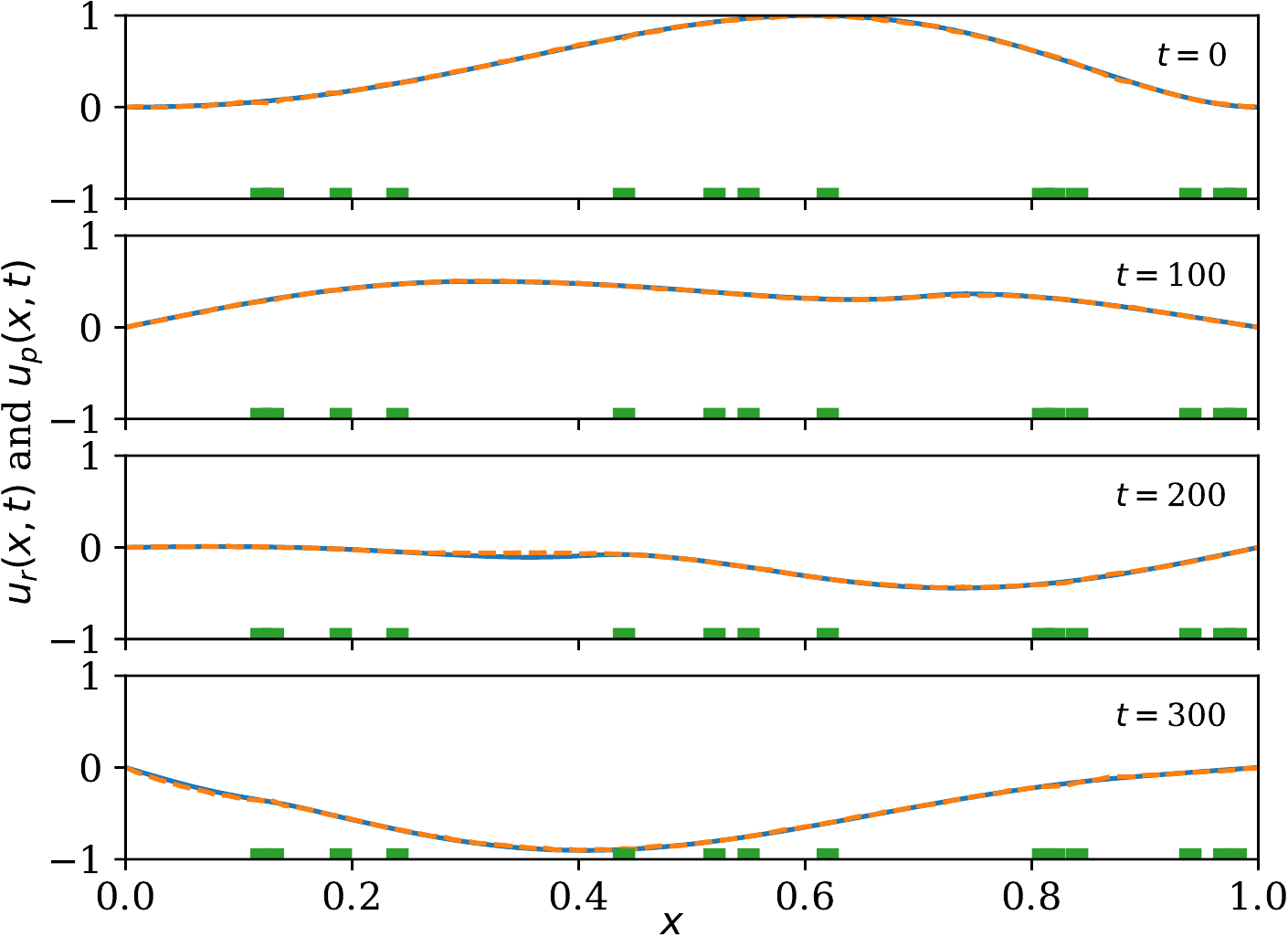}
\end{figure}

\begin{figure}[h!]
	\caption{Snapshots of prediction performance for the linear one-dimensional wave equation \eqref{1Dlinearwave} with initial condition $u_0(x) = x^2(x-1)^2(x-1/4)^2$, scaled to $[0,1]$. The solid blue line is the solution; the {orange} broken line is the prediction, and green squares show the optimized sensor arrangement.}
	\label{fig-snapshots-wave-IC3}
	\centering
	\includegraphics[width=0.5\textwidth]{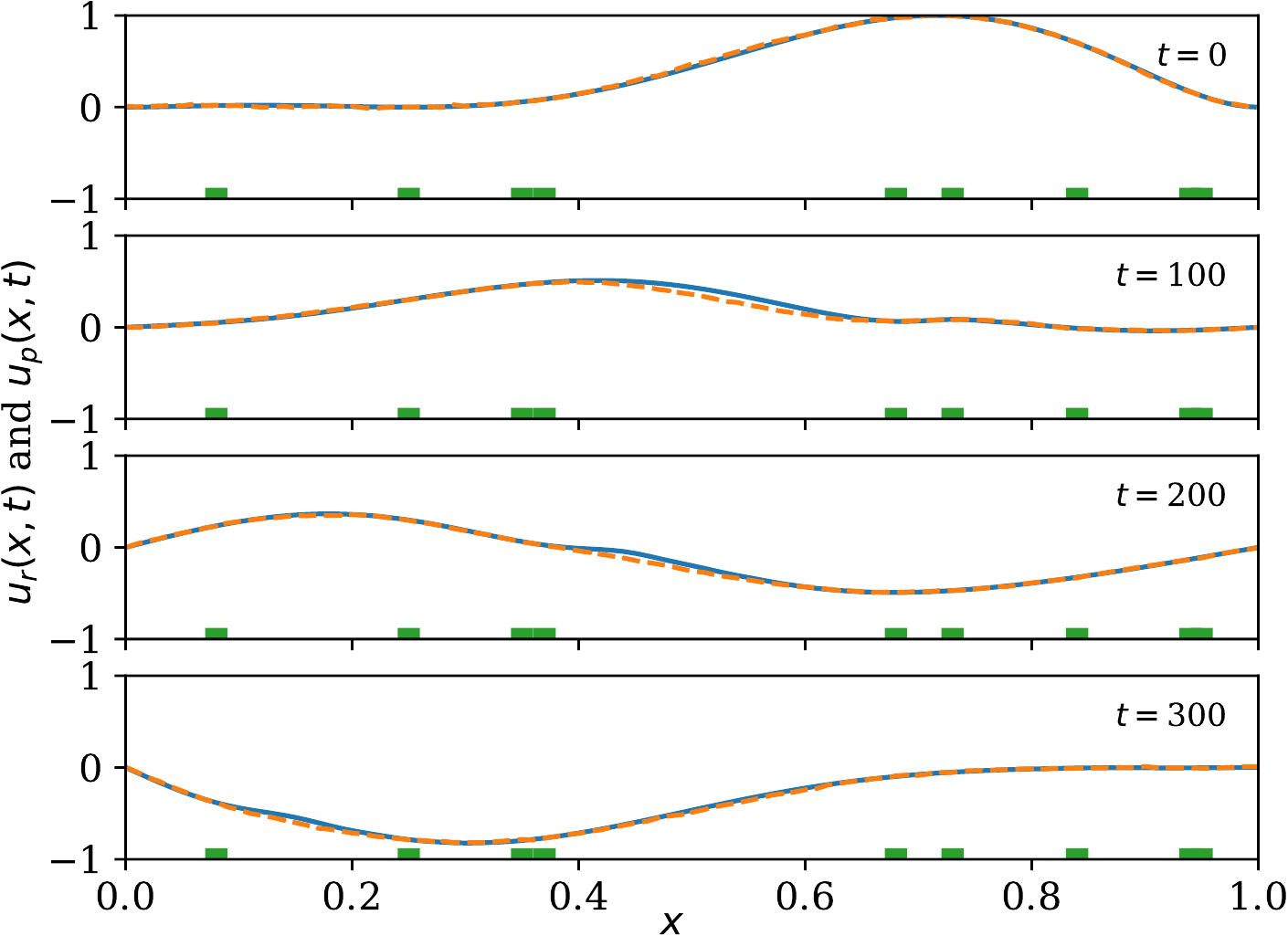}
\end{figure}

A time-centered and space-centered finite difference method is used to discretize the wave model and extract the solution data. 
The estimation performance is shown in \cref{fig-snapshots-wave-IC1,fig-snapshots-wave-IC2,fig-snapshots-wave-IC3} for various initial conditions $u_0$.
The sensor weight $\alpha$ is set to 2, {the mesh size $h$ is set to $10^{-2}$ (resulting in $c = 101$)}{; time increment $\Delta t$ is set to $0.1$;} and squared wave speed $\lambda$ is set to $3 \cdot 10^{-3}$.

\begin{figure}[h!]
	\caption{Reduction in cost for the linear one-dimensional wave equation \eqref{1Dlinearwave} and initial conditions given in \cref{fig-snapshots-wave-IC1,fig-snapshots-wave-IC2,fig-snapshots-wave-IC3}.}
	\label{fig-cost-wave-linear}
	\centering
	\includegraphics[width=0.5\textwidth]{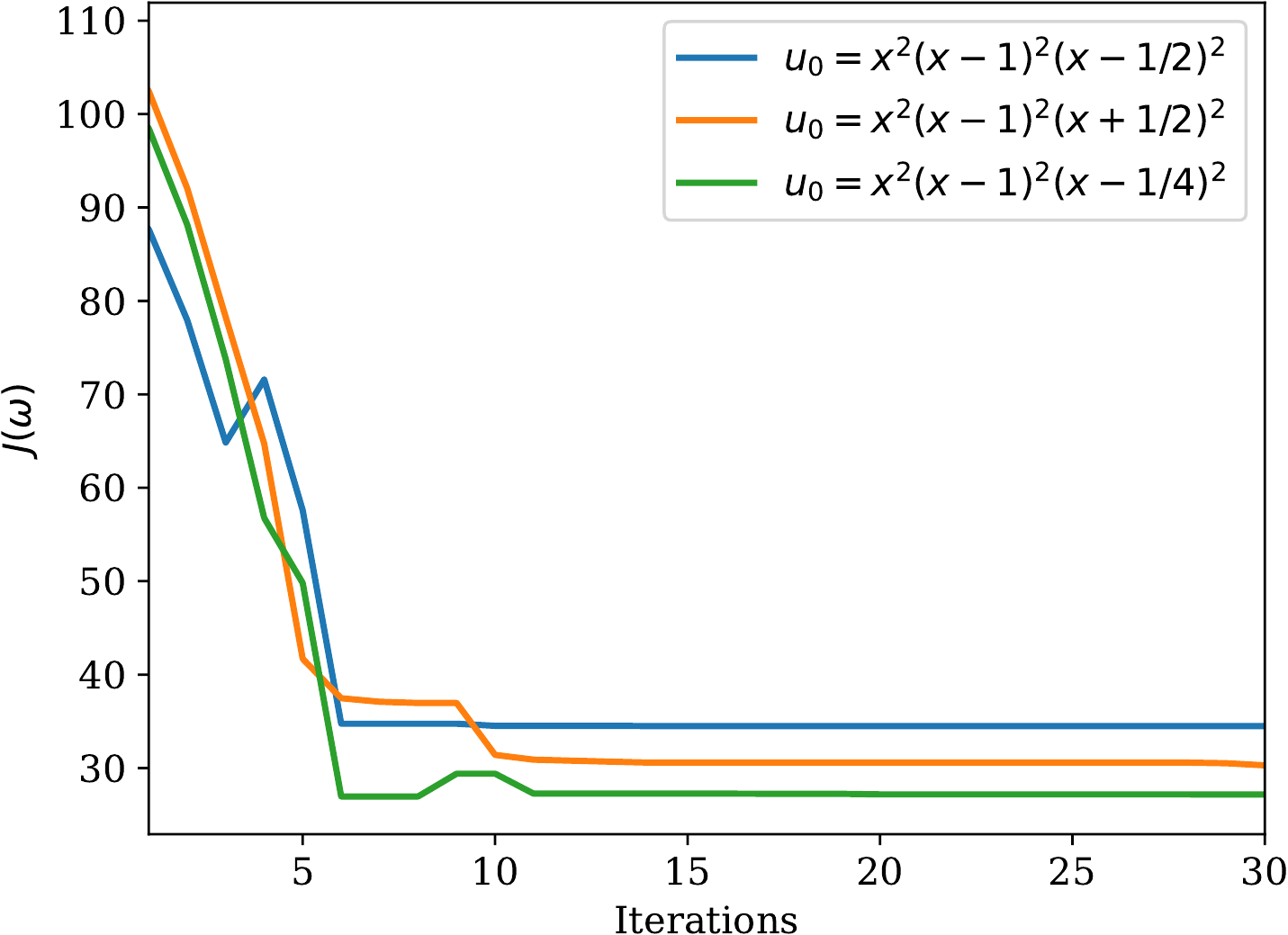}
\end{figure}

\begin{figure}[t!]
	\caption{$L^1$-norm of prediction error over time for the linear one-dimensional wave equation and all previous initial conditions. The norm of prediction error is larger than that of the heat equation since the solution to the heat equation converges to a steady state.}
	\label{fig-linftynorm-wave-linear}
	\centering
	\includegraphics[width=0.5\textwidth]{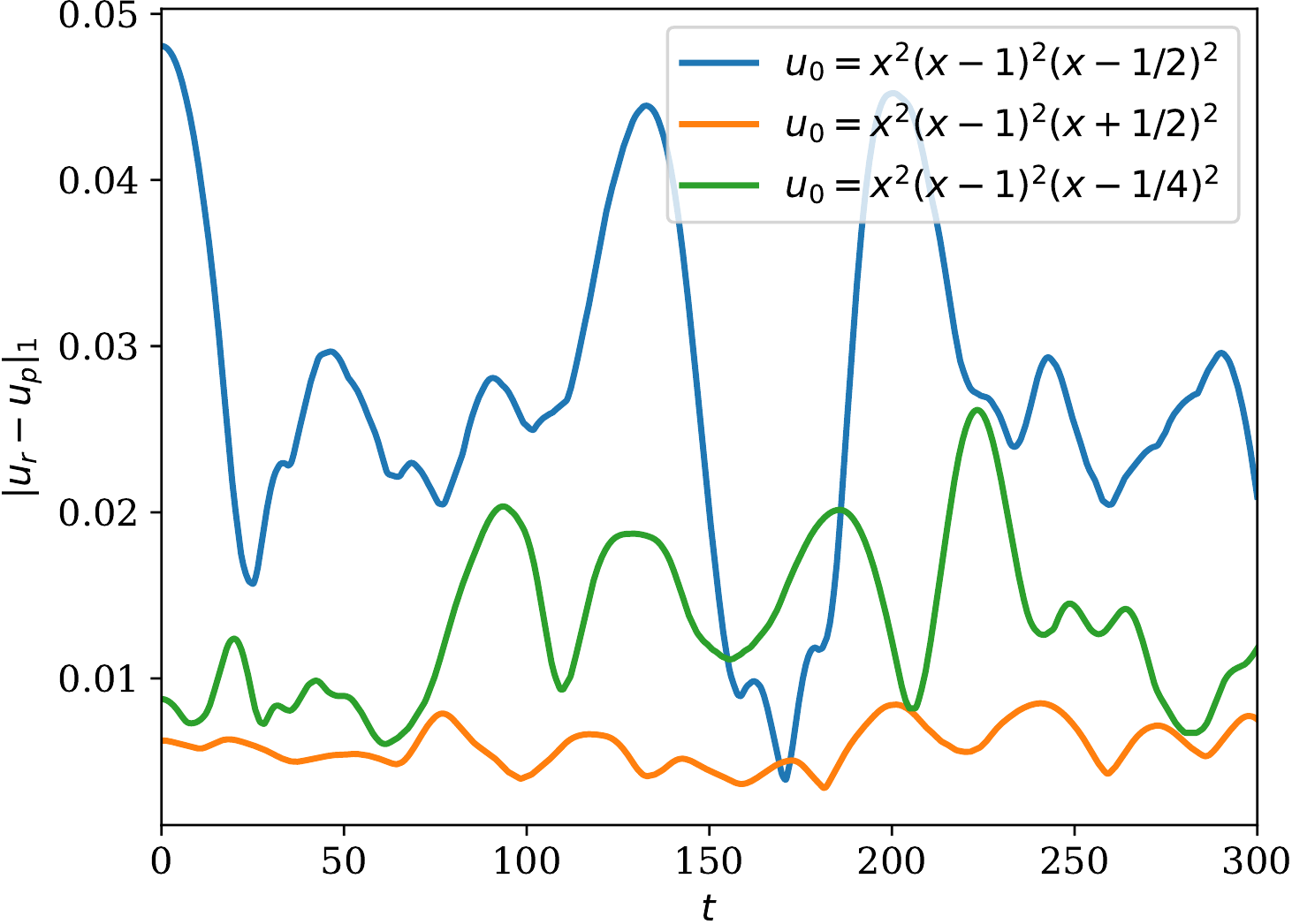}
\end{figure}

{Similar to the previous example, we observe good agreement between the predicted and the true state (\cref{fig-linftynorm-wave-linear}) although only very few sensors are in use.}
{In \cref{fig-snapshots-wave-IC1} and \cref{fig-snapshots-wave-IC3}, only 9\% of the wave region is covered with sensors, and in \cref{fig-snapshots-wave-IC2}, this number is 14\%.}
The reduction in cost over iteration is shown in \cref{fig-cost-wave-linear}.

\subsection{Two-dimensional heat equation}
\label{subsection:2Dlinearheat}
Consider the following heat equation in two space dimensions over $(x,y) \in [0,1]^2$
\begin{equation}\label{2Dlinearheat}
\begin{cases}
u_t (x,y,t) = \kappa \paren[auto](){u_{xx}(x,y,t) + u_{yy}(x,y,t)}, \\
u(x,y,0) = u_0,\\
u(0,y,t) = u_x (1,y,t) = 0,\\
u_y(x,0,t) = u_y(x,1,t) = 0.
\end{cases}
\end{equation}
The discretization is similar as in \cref{subsection:1Dlinearheat}.
The mesh size is uniform with length $h = 3\times 10^{-2}$; time increment $\Delta t$ is set to $0.1$; and conductivity $\kappa$ is set to $10^{-4}$. The initial condition is set to $u_0(x,y) = x^2y^2(x-1)^2(y-1)^2$, scaled to $[0,1]$.

\begin{figure}[h!]
	\caption{Optimized sensor arrangement for the two-dimensional linear heat equation \eqref{2Dlinearheat}. Green dots indicate the presence of a sensor.}
	\label{fig-optimal-sensor-2D-heat}
	\centering
	\includegraphics[width=0.5\textwidth]{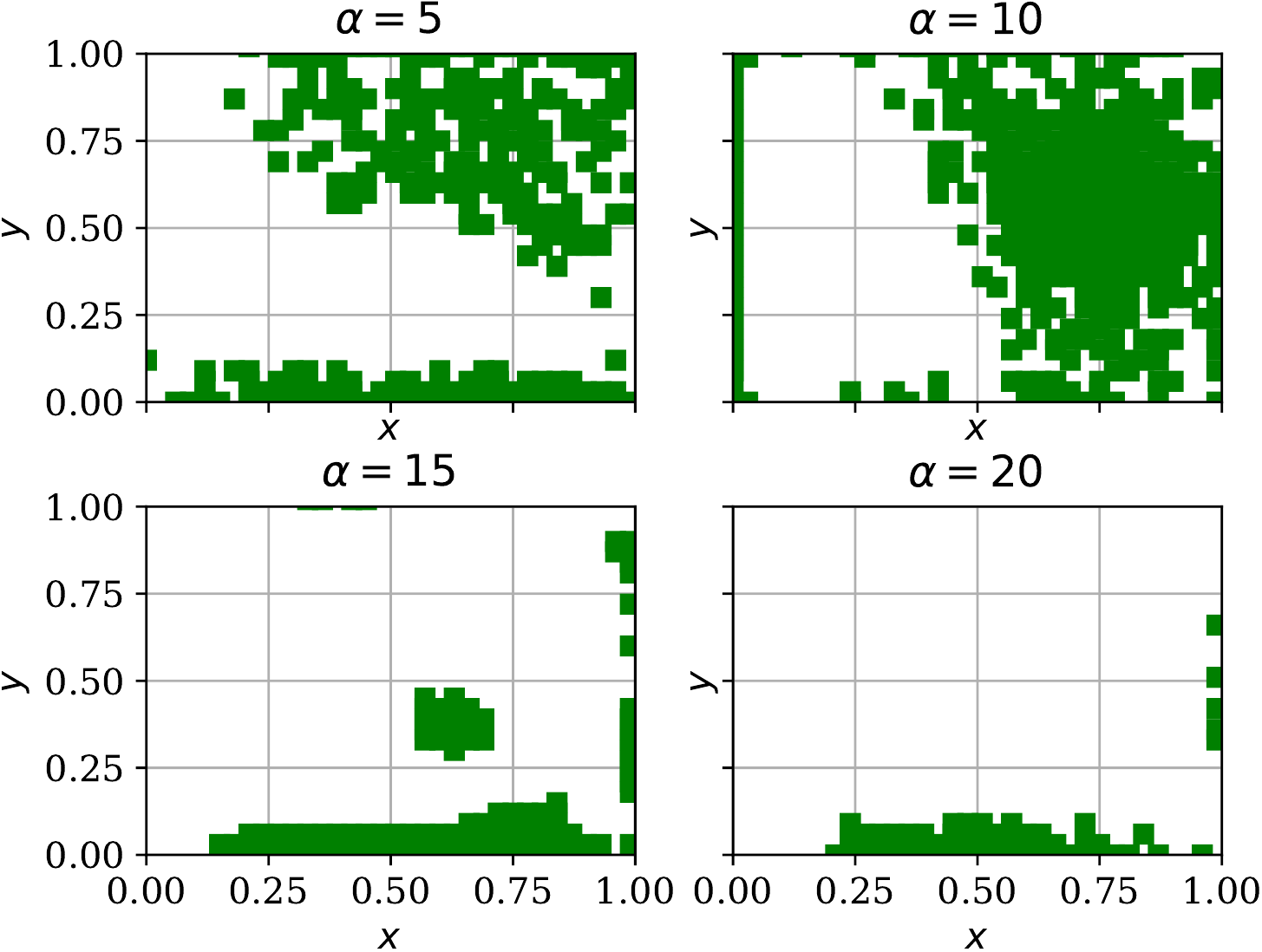}
\end{figure}
\vspace{0cm}
\begin{figure}[hbt!]
	\caption{Snapshots of prediction performance for two-dimensional linear heat equation \eqref{2Dlinearheat} and sensor weight $\alpha = 10$.}
	\label{fig-snapshots-2D-heat}
	\centering
	\begin{subfigure}[h!]{\textwidth}
		\includegraphics[width=0.5\textwidth]{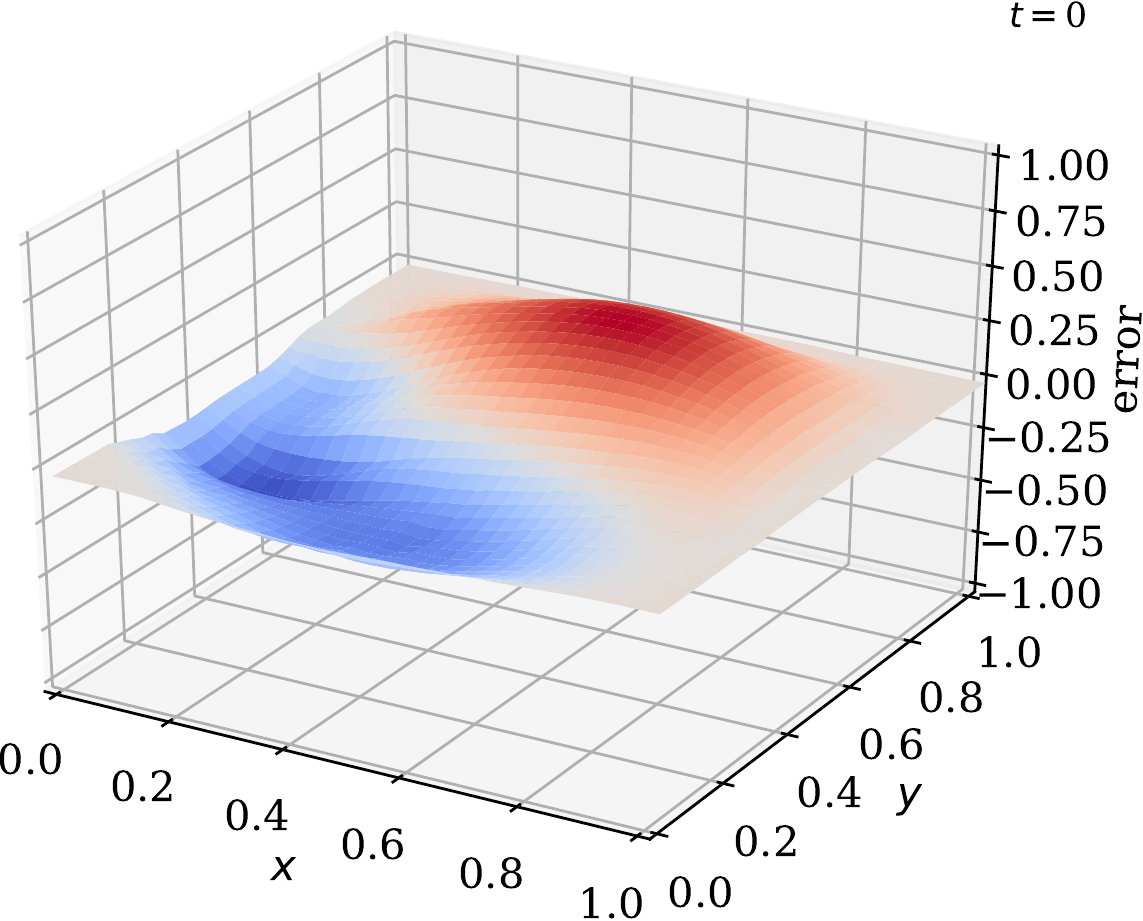}
	\end{subfigure}
	\begin{subfigure}[h!]{\textwidth}
		\includegraphics[width=0.5\textwidth]{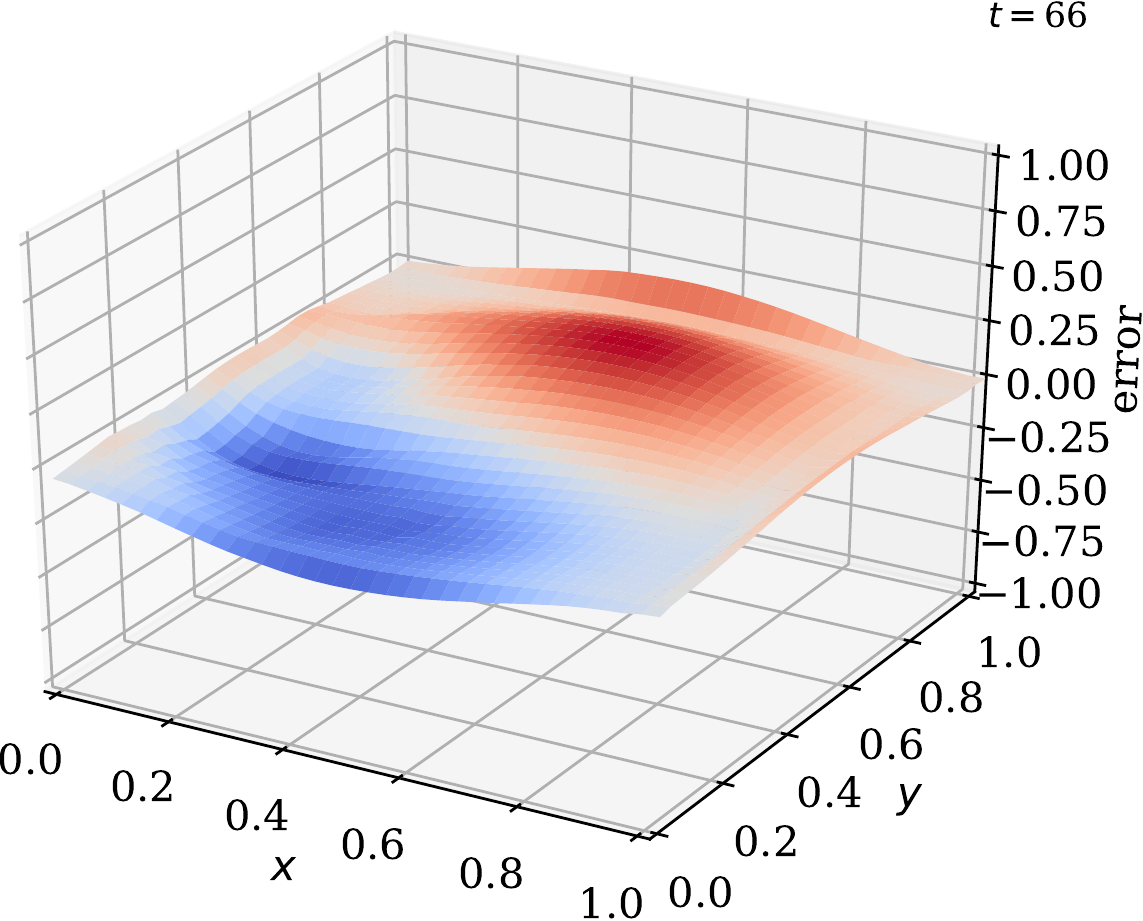}
	\end{subfigure}
	\begin{subfigure}[h!]{\textwidth}
		\includegraphics[width=0.5\textwidth]{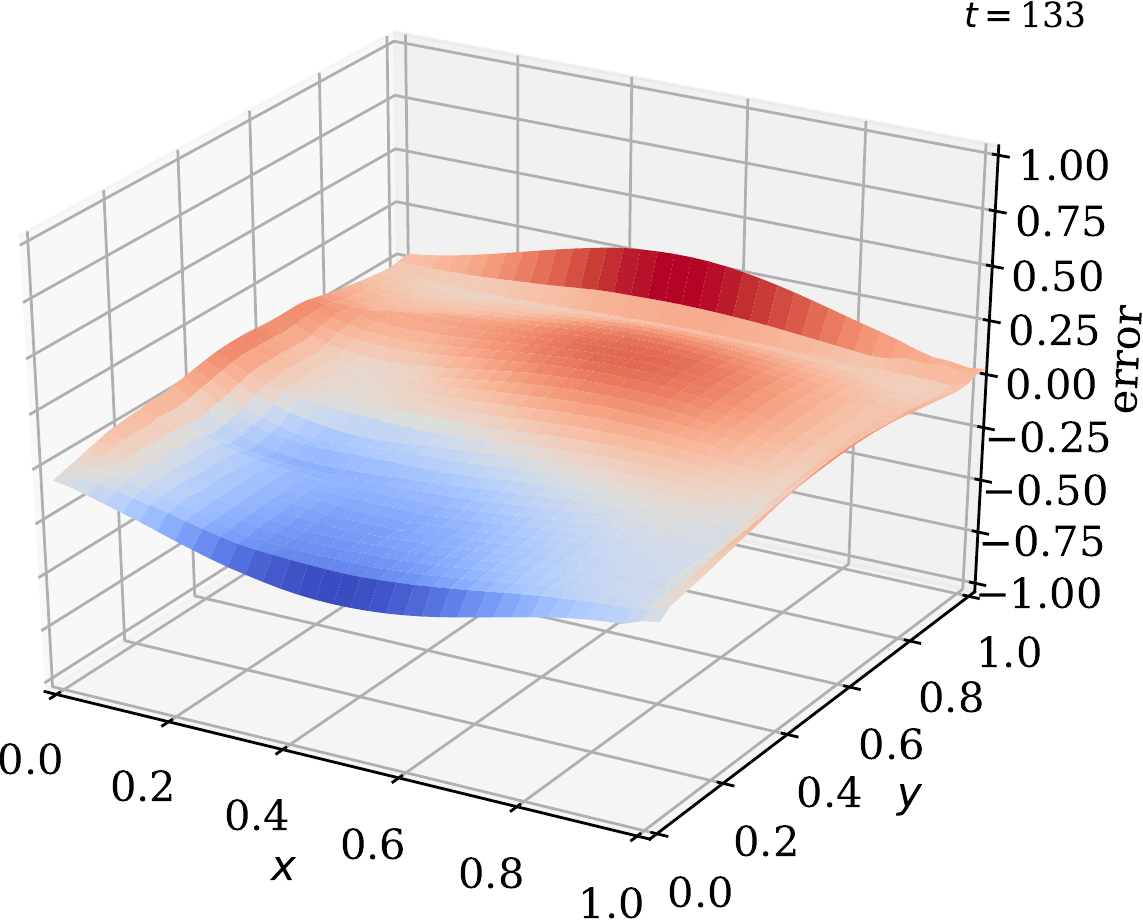}
	\end{subfigure}
\end{figure}

\begin{figure}[h!]\ContinuedFloat
	\begin{subfigure}[h!]{\textwidth}
		\includegraphics[width=0.5\textwidth]{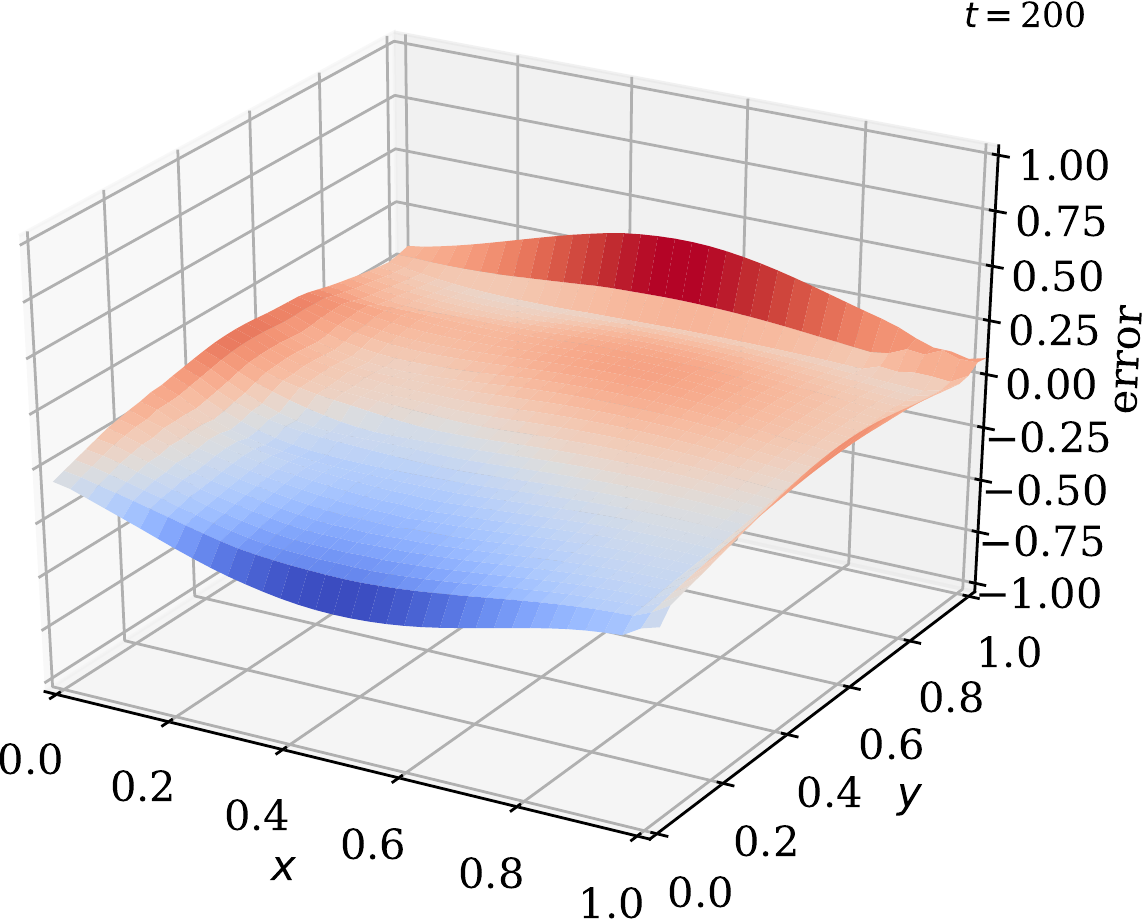}
	\end{subfigure}
\end{figure}

\begin{figure}[h!]
	\caption{Reduction in cost for the linear two-dimensional heat equation \eqref{2Dlinearheat} for different sensor weights.}
	\label{fig-cost-2D-heat}
	\centering
	\includegraphics[width=0.5\textwidth]{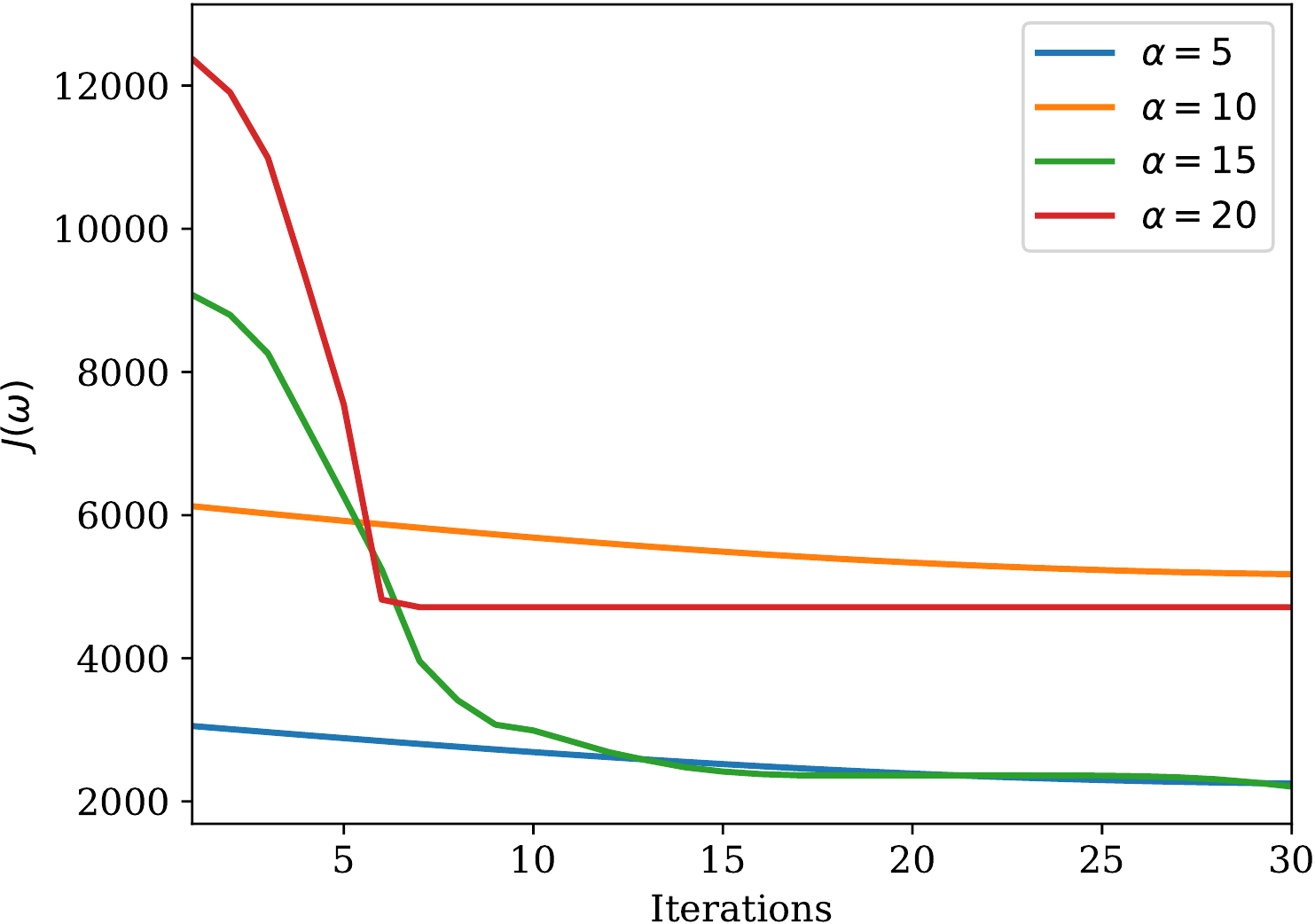}
\end{figure}

\begin{figure}[h!]
	\caption{$L^1$-norm of prediction error over time for the linear two-dimensional heat equation \eqref{2Dlinearheat} for different sensor weights.}
	\label{fig-norm-2D-heat}
	\centering
	\includegraphics[width=0.5\textwidth]{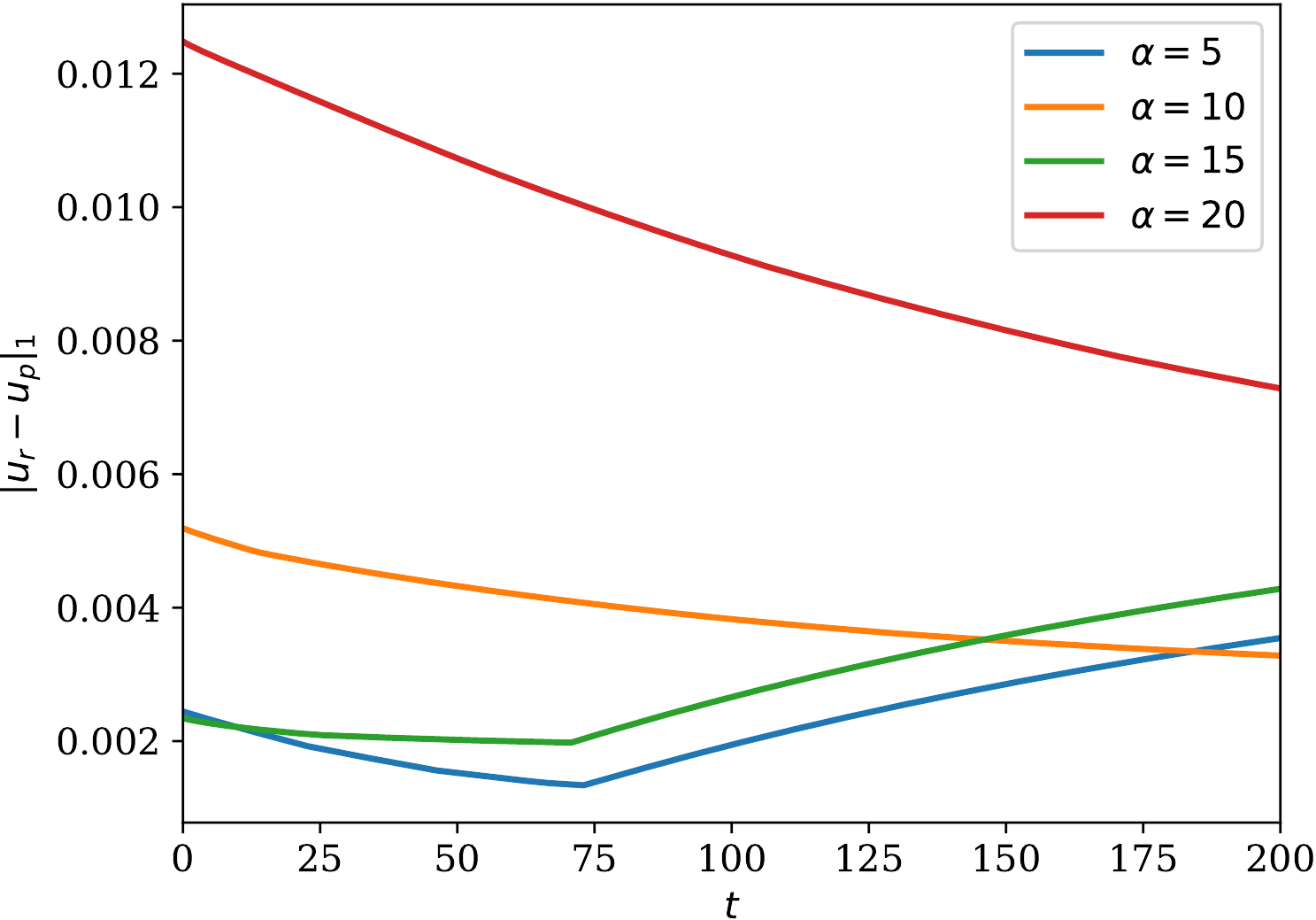}
\end{figure}

{In this example, we vary the value of the sensor weight coefficient $\alpha$.
	Once again, we observe good agreement between the predicted and the true state (\cref{fig-norm-2D-heat}) although only few sensors are in use; see \cref{fig-optimal-sensor-2D-heat,fig-snapshots-2D-heat}.}
The reduction in cost over iteration is shown in \cref{fig-cost-2D-heat}.

\section{Conclusion}
\label{section:conclusion}
Optimal prediction for distributed parameter systems was discussed in this paper. 
The optimal predictions involves two steps. 
In the first step, a neural network predictor is designed. 
In the second step, sensor shapes are optimized using the gradient of the network. 
Computer simulations are conducted for several models including heat and wave equation in one and two space dimensions. 
The results show a significant reduction in cost of prediction as well as improvement in the predictor performance. The computer codes are accessible in the repository \url{https://github.com/TUCMath/optimal-predictor}. 

\FloatBarrier

\ifdefined\isaccepted
\section*{Funding}
This work is part of a measure which is co-financed by tax revenue based on the budget approved by the members of the Saxon state parliament. Financial support is gratefully acknowledged.

\fi

\bibliographystyle{icml2021}
\bibliography{manuscript-icml2021.bib}

\end{document}